\documentclass[reqno]{amsart}

\usepackage{amsmath,amssymb,amsthm,amsfonts}
\usepackage{hyperref}
\usepackage{appendix}
\usepackage{graphicx}
\usepackage{setspace}
\usepackage{mathrsfs}
\usepackage{amsmath}

\usepackage{color}

\newtheorem{lemma}{Lemma}[section]
\newtheorem{theorem}{Theorem}[section]

\newtheorem{remark}{Remark}[section]

\numberwithin{equation}{section}

\begin{document}
\title[STABILITY OF INCOMPRESSIBLE PLANE COUETTE FLOW]{Stability for two-dimensional plane Couette flow to the incompressible Navier-Stokes equations with Navier boundary conditions}
\thanks{$^*$Corresponding author}
\author[Shijin Ding, Zhilin Lin]{Shijin Ding, Zhilin Lin$^*$}
\address[S. Ding]{South China Research Center for Applied Mathematics and Interdisciplinary Studies, South China Normal University,
Guangzhou, 510631, China}\address{School of Mathematical Sciences, South China Normal University,
Guangzhou, 510631, China}
\email{dingsj@scnu.edu.cn}
\address[Corresponding author: Z. Lin]{School of Mathematical Sciences, South China Normal University,
Guangzhou, 510631, China}
\email{zllin@m.scnu.edu.cn}

\date{\today}

\begin{abstract}
 This paper concerns with the stability of the plane Couette flow resulted from the motions of boundaries that the top boundary $\Sigma_1$ and the bottom one $\Sigma_0$ move with constant velocities $(a,0)$ and $(b,0)$, respectively. If one imposes Dirichlet boundary condition on the top boundary and Navier boundary condition on the bottom boundary with Navier coefficient $\alpha$, there always exists a plane Couette flow which is exponentially stable for nonnegative $\alpha$ and any positive viscosity $\mu$ and any $a, b \in \mathbb{R}$, or, for $\alpha<0$ but viscosity $\mu$ and the moving velocities of boundaries $(a,0), (b,0)$ satisfy some conditions stated in Theorem \ref{the1.1}. However, if we impose Navier boundary conditions on both boundaries with Navier coefficients $\alpha_0$ and $\alpha_1$, then it is proved that there also exists a plane Couette flow (including constant flow or trivial steady states) which is exponentially stable provided that any one of two conditions on $\alpha_0,\alpha_1$, $a, b$ and $\mu$ in Theorem \ref{the1.2} holds. Therefore, the known results for the stability of incompressible Couette flow to no-slip (Dirichlet) boundary value problems are extended to the Navier boundary value problems.
\end{abstract}

\maketitle

\vspace{-5mm}

\section{Introduction}\label{intro}

In this paper, we consider the stability of plane Couette flow for viscous incompressible fluid in a two dimensional slab domain, periodic in $x$ direction, $\Omega=\mathbb{T} \times (0,1) (\mathbb{T}=[-\pi,\pi])$ with the boundary $\Sigma=\Sigma_0 \cup \Sigma_1$, where $\Sigma_i=\{y=i\},\ i=0,1$. The motion of the incompressible fluid in $\Omega$ is governed  by the following incompressible Navier-Stokes equations
\begin{equation}\label{1.1}
\left\{
\begin{array}{ll}
\partial_t v-\mu \Delta v + v \cdot \nabla v +\nabla q=0,\\
\nabla \cdot v =0,
\end{array}
\right.
\end{equation}
where $v(t;x,y)=(v_1(t;x,y),v_2(t;x,y)) \in \mathbb{R}^2$ and $q(t;x,y)$ are the velocity and pressure, respectively. The constant $\mu>0$ is the viscosity.

To set our problem, we need to impose the boundary conditions. In this paper, we consider two cases, which are sated as follows.

{\bf Case I.} In the first case, the no-slip (Dirichlet) condition is imposed on the top boundary $\Sigma_1$ and the Navier condition is imposed on the bottom boundaries $\Sigma_0$. Since the Couette flow is resulted from the motion of the boundary, we suppose that the top boundary $\Sigma_1$ moves with a constant velocity $(a,0)$ and the bottom one $\Sigma_0$ with velocity $(b,0)$, where $a,b \in \mathbb{R}$ are constants (see \cite{6} for instance). That is,
\begin{equation}\label{1.2}
\left\{
\begin{array}{ll}
v \cdot \mathbf{n} =0 \ \mathrm{on} \ \Sigma,\\
v =(a,0) \ \mathrm{on} \ \Sigma_1,\\
\mathbb{S}(v) \cdot \mathbf{n} \cdot \mathbf{\tau} +\alpha (v-(b,0)) \cdot \mathbf{\tau}=0 \ \mathrm{on} \ \Sigma_0,
\end{array}
\right.
\end{equation}
where $\mathbb{S}(v)=-qI_2+\mu (\nabla v + \nabla^T v)$, $I_2$ is the $2\times 2$ identity matrix, $\mathbf{n}$ is the unit outward normal to the boundary and $\mathbf{\tau}$ is the tangential vector, and $\alpha$ is the constant of slip length. It should be pointed out that the term $v-(b,0)$ in condition (\ref{1.2}) represents the \emph{slip velocity}, see \cite{6} for more details.

It is well known that the Couette flow is an important type of shear flow in hydrodynamic stability theory.
In this case, it is direct to check that the Couette flow $(v_s, q_s)$ with
$$v_s=\left( \frac{\alpha (a-b)}{\mu +\alpha} y +\frac{\mu a+\alpha b}{\mu +\alpha},0 \right), \ q_s=\mathrm{constant}$$
is a steady solution to the problem (\ref{1.1})-(\ref{1.2}).

Let $u=v-v_s,\ p=q-q_s$. Then the Navier-Stokes equations around the Couette flow read as
\begin{equation}\label{1.3}
\left \{
\begin{array}{ll}
\partial_t u-\mu \Delta u + u \cdot \nabla v_s+v_s \cdot \nabla u +\nabla p=-u\cdot \nabla u \ \  & \mathrm{in} \ \Omega, \\
\nabla \cdot u =0 \ \ & \mathrm{in} \ \Omega,
\end{array}
\right.
\end{equation}
the corresponding boundary conditions are given as follows:
\begin{equation}\label{1.4}
\left \{
\begin{array}{lll}
u=0 \ & \mathrm{on} \ \Sigma_1,\\
u_2=0 \ & \mathrm{on} \ \Sigma_0, \\
\mu \partial_y u_1-\alpha u_1=0 \ & \mathrm{on} \ \Sigma_0.
\end{array}
\right.
\end{equation}

{\bf Case II.} In the second case, the Navier conditions are imposed on both boundaries, which are formulated as follows
\begin{equation}\label{1.5}
\left \{
\begin{array}{lll}
v \cdot \mathbf{n}=0 \ \ &\mathrm{on} \ \Sigma, \\
\mathbb{S}(v) \cdot \mathbf{n} \cdot \mathbf{\tau} + \alpha_1(v-(a,0)) \cdot \mathbf{\tau}=0 \ \ &\mathrm{on} \ \Sigma_1, \\
\mathbb{S}(v) \cdot \mathbf{n} \cdot \mathbf{\tau} + \alpha_0(v-(b,0)) \cdot \mathbf{\tau}=0 \ \ &\mathrm{on} \ \Sigma_0,
\end{array}
\right.
\end{equation}
in which $(a,0), (b,0)$ are the motion velocities of the boundaries $\Sigma_0, \Sigma_1$, respectively. The terms $v-(a,0), v-(b,0)$ are slip velocities.

In this case, for any $a, b \in \mathbb{R}$,
the problem, (\ref{1.1}) with (\ref{1.5}), admits a plane Couette flow $(v_s,q_s)$ with
$$v_s=\left(\frac{\alpha_0 \alpha_1(a-b)}{\mu(\alpha_0+\alpha_1) +\alpha_0\alpha_1} y +\frac{\mu(\alpha_1 a+\alpha_0 b)+\alpha_0\alpha_1 b}{\mu(\alpha_0 +\alpha_1)+\alpha_0\alpha_1},0\right)$$
and
$$q_s=\mathrm{constant}.$$

Therefore, we obtain the following perturbed problem
\begin{equation}\label{1.6}
\left \{
\begin{array}{ll}
\partial_t u +Lu=P(-u\cdot \nabla u) \ \  & \mathrm{in} \ \Omega, \\
u \cdot \mathbf{n}=0 \ \ & \mathrm{on} \ \Sigma, \\
\mathbb{S}(u) \cdot \mathbf{n} \cdot \mathbf{\tau} +\alpha_i u \cdot \mathbf{\tau}=0 \ \ &\mathrm{on} \ \Sigma_i, \ i=0,1,
\end{array}
\right.
\end{equation}
where $L$ is a linearized operator defined as
$$Lu=P(-\mu \Delta u + u \cdot \nabla v_s+v_s \cdot \nabla u),$$
and $P$ is the Helmholtz projection (see Section \ref{sec2}).

\begin{figure}[htbp]\label{1.8}
  \centering
  \includegraphics[width=0.66\textwidth]{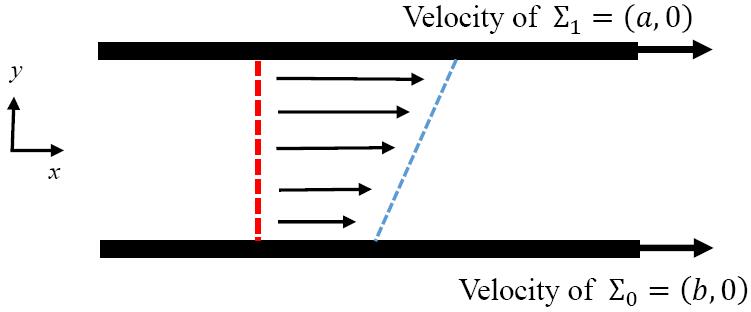}
  \caption{\rm Couette flow in $\mathbb{R} \times (0,1)$ ($a,b>0$)}\label{fig:digit}
\end{figure}

The Couette flows with $a, b>0$ can be shown in Figure \ref{1.8} (Note that $a, b \in \mathbb{R}$). It should be pointed out that in Case I-II, $v_s$ is reduced to the constant flow $v_s=(a,0)$ provided that $a=b$.

Our aim is to study the asymptotic stability for the nonlinear problem (\ref{1.3})-(\ref{1.4}) and (\ref{1.6}).

Let us give a brief review about the stability theory and some related problems.
The stability of trivial or non trivial steady states to the equations (\ref{1.1}) with no-slip (Dirichlet) boundary conditions have been studied for a long time. For the stability of trivial steady states such as Rayleigh-Taylor stability and instability, we refer to \cite{YGuo,Jiang1,Jiang2,Wang}. However, the researches for the stability of non trivial steady states such as Couette flows, Poiseuille flows or general shear flows are far from completion. The first result for the stability of incompressible plane Couette flows with no-slip (Dirichlet) boundary conditions was obtained by Romanov in a beautiful paper \cite{Romanov}, which shows that the plane Couette flow is stable for any fixed Reynolds number. Similar result was obtained by Heck et al. \cite{Heck} for periodic case. For the general shear flows including Poiseuille flow, in the large Reynolds number regime, the spectral instability was obtained by Grenier, Guo and Nguyen \cite{Guo}. The nonlinear stability for the cylindrically symmetric Poiseuille flow was obtained by Gong and Guo \cite{Gong}.

It is well known that the key point for stability and instability problems is the spectral analysis. For the linearized operator around the trivial steady states, the spectral analysis can be done by searching for growing normal mode solutions and using variational method (see \cite{YGuo,Jiang1,Jiang2,Wang}) since the linearized operators are self-adjoint. However, the linearized operators around the shear flows are always non self-adjoint and nonlocal. The spectral analysis mainly depends on the analysis of the Orr-Sommerfeld equations. For any steady shear flow $v_s=U(y)e_1$, by using the normal Fourier transform $\psi=\phi(y) e^{ik (x-ct)}, k \in \mathbb{R}, c \in \mathbb{C}$, where $\psi$ is the stream function such that $u=\nabla^\perp \psi,$ one can get the Orr-Sommerfeld equation around the shear flow $v_s$
\begin{equation}\label{1.9}
(\partial_y^2-k^2)^2 \phi=ikR\left[\left(U-c\right)(\partial_y^2-k^2) \phi-U''\phi\right]
\end{equation}
with suitable boundary conditions, where $R$ is the Reynolds number. By the classical spectral stability theory, the flow $v_s$ is linearly spectral stable for $\mathrm{Im}\ c<0$ for any $c\in \mathbb{C}$ and unstable for $\mathrm{Im}\ c_0>0$ for some $c_0\in\mathbb{C}$.

The study for the Orr-Sommerfeld equation was initialed by Orr in 1907, see \cite{Orr1,Orr2} for details. Up to now, there are few results about Orr-Sommerfeld equations with no-slip boundary conditions, see \cite{Drazin,Lin1,Lin2} for instance. For the spectrum analysis of the Orr-Sommerfeld equation, Joseph \cite{Joseph1,Joseph2} gave the eigenvalue bounds for the Orr-Sommerfeld equation, which established some sufficient conditions for stability. Some other similar results can be found in \cite{Yih,Adelina}.

For the stability problems in compressible Navier-Stokes equations with no-slip boundary conditions, most results are also obtained via the spectral analysis of the linearized perturbation operator. A sufficient condition for the stability of the compressible Couette flow was obtained by Kagei \cite{Kagei1}. With the similar idea, Kagei and Nishida \cite{Kagei2} proved that the Poiseuille flow is unstable if Reynolds number and Mach number satisfy some conditions. Recently, Li and Zhang \cite{HLLi} improved the result of \cite{Kagei1}.

It is interesting to compare Navier boundary conditions with the no-slip boundary conditions in our problem. The no-slip (Dirichlet) boundary conditions mean that the fluid does not slip along the boundary.
However, this is not always realistic and leads to the strong boundary layer in general. For example, hurricanes and tornadoes, do slip along the ground, lose energy as they slip and do not penetrate the ground. Other examples about the slip of the fluid on the boundary occur when moderate pressure is involved such as in high altitude aerodynamics, or in immiscible two phase flows, the moving contact line is not compatible with no-slip boundary condition. To describe these phenomenons, Navier\cite{Navier} in 1823 introduced the so called Navier boundary conditions. The Navier boundary condition is formulated as
$$v \cdot \mathbf{n} =0,\ \ \mathbb{S}(v) \cdot \mathbf{n} \cdot \mathbf{\tau} +\alpha v \cdot \mathbf{\tau}=0 \ \mathrm{on} \ \partial{\Omega},$$
in which $\alpha$ is a physical parameter standing for the frictions between the fluid and the ground or permeability and others which is either a constant
or a $L^{\infty}(\partial{\Omega})$ function\cite{Kelliher}, even a smooth matrix\cite{Gie}.

The case $\alpha \geq 0$ is the classical case which reflects the friction between the fluid and the boundary and has got extensive attentions by physicists and mathematicians in studying the existence, uniqueness,
regularity and vanishing viscosity to system (\ref{1.1}), see for instance \cite{Xin1,Xin2}. However, the case $\alpha<0$ does exist in reality and in physics. For example, for flat hybrid
gas-liquid surfaces, the effective slip length $\alpha$ is always negtive\cite{Haase}. Navier boundary condition with $\alpha<0$ is also used for the simulations of flows in the presence of rough boundaries such as in aerodynamics, or in the case of permeable boundary in which the Navier boundary condition was called Beavers-Joseph's law \cite{Amrouche,6}, or in weather forecasts and in hemodynamics \cite{6,1993}, or when the boundary wall accelerates the fluid \cite{3,1995}.

In this paper, we assume that $\alpha, \alpha_0$ and $\alpha_1$ are constants.

J.-L. Lions\cite{JLions} and P.-L. Lions\cite{PLions} considered the following boundary conditions, which is called vorticity free boundary conditions:
$$v \cdot \mathbf{n} =0,\ \ \omega(v)=0 \ \mathrm{on} \ \partial{\Omega},$$
where $\omega(v)=\partial_x v_1-\partial_y v_2$ is the vorticity of $v$. In other words, the vorticity free boundary condition is the special case of the Navier boundary condition when $\frac{\alpha}{\mu}=2\kappa$, where $\kappa$ is the curvature of the boundary $\partial{\Omega}$, see for instance \cite{JLions,Kelliher}. Therefore, for our problem, Navier boundary conditions contain vorticity free boundary condition provided that $\alpha=0$ or $\alpha_0=\alpha_1=0$.

In view of the results of Romanov \cite{Romanov} and Heck \cite{Heck}, it is very natural to consider the stability problem with Navier boundary conditions. In our results, for the Navier boundary problem, we can find some sufficient conditions for the stability of Couette flow. The sufficient conditions depend on the viscosity $\mu$, the moving velocities of boundaries $(a,0), (b,0)$ and the Navier coefficients $\alpha$ or $\alpha_0,\ \alpha_1$.

Similar results for the stability and instability of trivial steady states $(0,q_s)(q_s=\mathrm{constant})$ with Navier boundary conditions were obtained by the first author, Li and Xin \cite{SJDing}, which provided a critical viscosity determined by the Navier coefficients to distinguish the stability from the instability. In addition, in \cite{SJDing}, the Navier boundary condition with $\alpha\geq 0$ is called {\it dissipative} and the Navier boundary condition with $\alpha<0$ is called {\it absorptive}.

Our aim is to analyze the stability of the incompressible Couette flow with Navier boundary conditions. One key step is to determine the sign of the image part of spectrum for the Orr-Sommerfeld equation. The key point is to estimate the upper bound of $\mathrm{Im} \ c$. Therefore, we need to establish estimates for Orr-Sommerfeld equation. Compared with the cases in Joesph \cite{Joseph1,Joseph2} and Romanov \cite{Romanov}, we have to deal with the boundary terms resulted from the Navier boundary conditions. To overcome the difficulties, we will modify the idea of Joseph \cite{Joseph1,Joseph2} and obtain the desired estimates.

For {\bf Case I}, if $\alpha \geq 0,$ our main results imply that the Couette flow is asymptotically nonlinear stable under small perturbation for any viscosity $\mu >0$ and any moving velocities of the boundaries $(a,0)$ and $(b,0)(\forall a,b \in \mathbb{R}).$ That is, the results of Romanov\cite{Romanov} still hold if $\alpha \geq 0.$ However, if $\alpha < 0$, our main results yield that the Couette flow is asymptotically nonlinear stable for small perturbation provided that $\alpha$ and $\mu$ satisfy the conditions that $\mu>-3\alpha$ and $\frac{|\alpha(a-b)|}{\mu(\mu+\alpha)}(1+\frac{3\alpha}{\mu})^{-\frac{1}{2}}<2\sqrt{2},$ see Theorem \ref{the1.1}. In addition, this result implies that the Couette flow is stable for all positive viscosity with vorticity free boundary conditions, see Remark \ref{rem1}.

 For {\bf Case II}, we can give a sufficient condition for stability, see Theorem \ref{the1.2}. If $\alpha_0,\alpha_1 \geq 0$, we show that the steady flow is asymptotically nonlinear stable under small perturbation for any viscosity $\mu >0$ and $a,b \in \mathbb{R}$. Therefore the results of Romanov\cite{Romanov} still hold if $\alpha_0 \geq 0$ and $\alpha_1 \geq 0$. Otherwise, the Couette flow is asymptotically nonlinear stable under small perturbation provided that $\alpha_0,\alpha_1, \mu$ and $a, b$ satisfy some conditions, see Theorem \ref{the1.2}.

In the {\bf Case II}, it should be noted that Couette flows is reduced to the trivial steady state $(v_s,q_s)=(0,\mathrm{constant})$ provided that $a=b=0$ or $\alpha_0=\alpha_1=0$, and the case of $a=b=0$ was studied by the first author, Li and Xin \cite{SJDing} recently. If $a=b=0$, for the trivial steady state $(v_s,q_s)=(0,\mathrm{constant})$, if $\alpha_0 \geq 0$ and $\alpha_1\geq 0$, then the Theorem \ref{the1.2} implies that the steady state is stable for any viscosity $\mu>0$, which is the same as in \cite{SJDing}. Otherwise, the steady state $(v_s,q_s)=(0,\mathrm{constant})$ is stable provided that the condition (iii) of Theorem \ref{the1.2} holds (the condition (iv) holds surely since $a=b=0$). In addition, the first author, Li and Xin \cite{SJDing} gave a critical viscosity $\mu_c$ and they proved that the steady state $(v_s,q_s)=(0,\mathrm{constant})$ is stable provided that $\mu >\mu_c$. Here we can not obtain such a critical viscosity to distinguish the stability from instability.

To state our results, let us introduce some notions and function spaces. The domain symbol $\Omega$ will be omitted for simplicity. Let
$$\mathscr{D}:=\left\{u(x,y)=\sum_{k \in J} \hat{u}_k(y) e^{ikx}: J \subset \mathbb{Z} \ \mathrm{is \ some \ finite \ subset}, \ \hat{u}_k(y) \in  C^{\infty}([0,1]) \right\}$$
and
$$\mathscr{D}_\sigma:=\left\{u\in \mathscr{D}: \nabla \cdot u=0 \right\},$$
where
$$\hat{u}_k(y)=\frac{1}{2\pi}\int_{\mathbb{T}} u(x,y) e^{-ikx}\mathrm{d}x, \ k \in \mathbb{Z}.$$

For the boundary conditions (\ref{1.4}), we define
$$\mathscr{D}_{*}:=\left\{ u \in \mathscr{D}:u \ \mathrm{satisfies \ the \ boundary \ conditions} \ (\ref{1.4})\right\}$$
and
$$\mathscr{D}_{*,\sigma}:=\left\{ u \in \mathscr{D}_\sigma: u \ \mathrm{satisfies \ the \ boundary \ conditions} \ (\ref{1.4})\right\}.$$

Define the norms
$$\Vert u \Vert_{L^p}^p=\int_{0}^{1} \int_{\mathbb{T}} |u(x,y)|^p \mathrm{d}x\mathrm{d}y$$
and
$$\Vert u \Vert_{W^{m,p}}^p=\sum_{|l| \leq m}\Vert D^l u \Vert_{L^p}^p. $$

With the above definitions, we can define the Sobolev spaces as the closures of $\mathscr{D},\mathscr{D}_{0,\sigma}$ or $\mathscr{D}_{*}$ with the following norms:
$$W^{m,p}=\overline{\mathscr{D}}^{\Vert \cdot \Vert_{W^{m,p}}}, \ L^p=W^{0,p},\ L^p_\sigma=\overline{\mathscr{D}_\sigma}^{\Vert \cdot \Vert_{L^p}},\ W^{m,p}_*=\overline{\mathscr{D}_*}^{\Vert \cdot \Vert_{W^{m,p}}},\ W^{m,p}_{*,\sigma}=\overline{\mathscr{D}_{*,\sigma}}^{\Vert \cdot \Vert_{W^{m,p}}},$$
and we denote
$$H^m=W^{m,2},\ H^m_*=W^{m,2}_*, \ H^m_{*,\sigma}=W^{m,2}_{*,\sigma}$$
for simplicity.

For the operator $L$, denote the spectrum of $-L$ by $\sigma(-L)$ and the resolvent set of $-L$ by $\rho(-L)$. In addition, for any $\theta >0$, define the sector of angle $\theta$ as
$$\Sigma(\theta):=\left\{ z \in \mathbb{C}\setminus\{0\}: |\mathrm{arg} \ z |<\theta\right\}.$$

For the problem (\ref{1.3})-(\ref{1.4}), our main result reads as follows.
\begin{theorem}\label{the1.1}
The Couette flow $v_s=\left( \frac{\alpha (a-b)}{\mu +\alpha} y +\frac{\mu a+\alpha b}{\mu +\alpha},0 \right)$ is linearly stable provided that any one of the followings (i), (ii) holds:

(i) $\alpha \geq 0, a,b \in \mathbb{R}$ and $\mu >0$;

(ii) $\mu>-3\alpha>0(i.e.,\alpha<0)$ and $\frac{|\alpha(a-b)|}{\mu(\mu+\alpha)}\cdot(1+\frac{3\alpha}{\mu})^{-\frac{1}{2}}<2\sqrt{2}.$

In addition, there exists $\varepsilon >0$ small enough such that if the initial data $ u_0 \in H^1_{*,\sigma}$ and $\left\Vert u_0 \right\Vert_{H^1} \leq \varepsilon$, then the problem (\ref{1.3})-(\ref{1.4}) is nonlinearly stable, i.e., there exist unique global solution $(u,p) \in (H^1_{*,\sigma} \cap H^2) \times H^1$ satisfy (\ref{1.3})--(\ref{1.4}), and the following decay holds
\begin{equation}\label{1.11}
\left\Vert u(t) \right\Vert_{H^1} \leq C_1e^{-\beta t} \left\Vert u_0 \right\Vert_{H^1},
\end{equation}
where the positive constants $C_1, \beta$ depend only on $\mu, \alpha, a, b$.
\end{theorem}

\begin{remark}\label{rem1}
Theorem \ref{the1.1} implies that the results of Romanov \cite{Romanov} still hold for the Navier boundary condition if $\alpha \geq 0$. In particular, let $\alpha=0$, then the Couette flow is reduced to a constant flow and the Navier boundary conditions become into vorticity free boundary conditions. In this case, of course, the results of Romanov \cite{Romanov} also hold for the vorticity free boundary conditions.
\end{remark}

For the problem (\ref{1.6}), we have the following result.
\begin{theorem}\label{the1.2}
The Couette flow $v_s=\left(\frac{\alpha_0\alpha_1(a-b)}{\mu(\alpha_0+\alpha_1) +\alpha_0\alpha_1} y +\frac{\mu(\alpha_1 a+\alpha_0 b)+\alpha_0\alpha_1 b}{\mu(\alpha_0 +\alpha_1)+\alpha_0\alpha_1},0\right)$ is linearly stable provided that any one of the followings (iii), (iv) holds:

(iii) $\alpha_0 \geq 0,\alpha_1\geq 0, a,b \in \mathbb{R}$ and $\mu>0$;

(iv) otherwise,

$\mu>\max\left\{(1+C_P)\max\limits_{l=0,1}\{|\alpha_l|\}-C_P(\alpha_0+\alpha_1),
2\max\limits_{l=0,1}\{|\alpha_l|\}-(\alpha_0+\alpha_1)\right\}$

\noindent and
$$\left|\frac{\alpha_0 \alpha_1(a-b)}{\mu\left(\mu(\alpha_0+\alpha_1) +\alpha_0\alpha_1\right)}\right|\cdot\left(1-\frac{2\max\limits_{l=0,1}|\alpha_l|-\alpha_0-\alpha_1}{\mu}\right)^{-\frac{1}{2}} <2\sqrt{2},$$

\noindent where constant $C_P>0$ is the best constant so that Poincar$\acute{\mathrm{e}}$ inequality for $f(y)\in H^1(0,1)$ with $\int_0^1f(y)\mathrm{d}y=0$ holds, see Lemma \ref{lemma5.1}.

In addition, there exists $\varepsilon >0$ small enough such that if $\left\Vert u_0 \right\Vert_{H^1} \leq \varepsilon,\ u_0 \in H^1_{**,\sigma}$, then the Couette flow is nonlinearly stable, i.e., there exists unique global solution $(u,p) \in (H^1_{**,\sigma} \cap H^2) \times H^1$ to (\ref{1.6}), and the following decay estimate holds
\begin{equation}\label{1.12}
\left\Vert u(t) \right\Vert_{H^1} \leq C_2e^{-\gamma t} \left\Vert u_0 \right\Vert_{H^1},
\end{equation}
where the positive constants $C_2, \gamma$ depend only on $\mu, a, b, \alpha_l(l=0,1), \Omega$ and $H^1_{**,\sigma}$ is defined in Section \ref{sec5}.
\end{theorem}

\begin{remark}\label{rem2.01}
Theorem \ref{the1.2} shows that the results of Romanov\cite{Romanov} also hold for the Navier boundary value problems if $\alpha_0 \geq 0,\alpha_1\geq 0$. Of course, this case includes the trivial steady state and the constant flow (Note that $v_s$ is reduced to constant flow if $a=b$).

In particular, if $\alpha_0=\alpha_1=0$, then $v_s=(0,0)$ is a trivial steady state. In this case, the results of Romanov\cite{Romanov} hold for the vorticity free boundary conditions, which is similar to Remark \ref{rem1}.
\end{remark}

Both above theorems give some sufficient conditions for the stability of the Couette flow in two cases. As mentioned before, the Couette flow is resulted from the motion of boundary, therefore together with the viscosity and slip length, the velocity of the motion should be concerned as the factor for stability or instability. Precisely, the \emph{relative velocity} $(a-b,0)$, the difference of motion velocities of two boundaries, will effect the energy of fluids with viscosity and slip length. According to our results, if the Navier boundary conditions are \emph{dissipative}, that is $\alpha \geq 0$ or $\alpha_0, \alpha_1 \geq 0$, then any motion velocities of the boundaries can not result in the instability, which means that the effect of slip lengths will be treated as the main factor for the stablity.
However, if the Navier boundary conditions are \emph{absorptive}, i.e., $\alpha <0$ or at least one of $\alpha_l(l=0,1)<0$, the stability of fluid will mainly depend on the viscosity. In other words, the viscosity should not be too small, or the modulus of relative velocity $|a-b|$ should not be too large.

The rest of this paper is organized as follows. In Section \ref{sec2}, we will introduce some elementary conclusions and inequalities which will be used in later analysis. Section \ref{sec3} is devoted to the proof of linear stability in Theorem \ref{the1.1}. The nonlinear stability in Theorem \ref{the1.1}
is shown in Section \ref{sec4}. In Section \ref{sec5}, we will prove the Theorem \ref{the1.2}.

\section{Preliminary}\label{sec2}

To define the Stokes operator and the perturbed operator $L$, we need some results about the Helmholtz projection and the resolvent problem, which ensure that the perturbed operator is well-defined and generates an analytic semigroup. The results can be obtained by applying the classical Fourier analysis, see \cite{Abel1,Abel2} for instance.

\begin{lemma}{\bf(\cite{Heck})}\label{lemma2.1}
For any vector field $u \in L^2,$ there exists unique vector field $v \in L^2_\sigma,$ such that
\begin{equation}\label{2.2}
u=v+\nabla p
\end{equation}
for some scalar $p \in H^1$. In addition, the following estimate holds
\begin{equation}\label{2.3}
\left\Vert v \right\Vert_{L^2} +\left\Vert \nabla p \right\Vert_{L^2} \leq C\left\Vert u \right\Vert_{L^2},
\end{equation}
where the constant $C>0$ depends only on $\Omega$.
\end{lemma}

\begin{remark}
The Lemma \ref{lemma2.1} implies that the Helmholtz projection
$$P: u \in  L^2 \mapsto v=Pu \in L^2_\sigma$$
is a bounded linear operator.
\end{remark}

By the Helmholtz projection, we define the Stokes operator $-A$ in $L^2_\sigma$ by
$$Au=P(-\mu \Delta u),\ \ \ u \in D(A)=H^1_{*,\sigma} \cap H^2.$$
Obviously, the operator $-A$ is unbounded in $L^2_\sigma$. And the following resolvent result for $-A$ is important.

\begin{lemma}\label{lemma2.3}
Suppose that $\theta\in (0,\frac{\pi}{2})$ and $\lambda \in \Sigma(\frac{\pi}{2}+\theta)$. Then for any $f \in L^2_\sigma$, there exists unique $u \in D(A)$ such that
\begin{equation}\label{2.4}
(\lambda+A)u=f,
\end{equation}
and the following estimate holds
\begin{equation}\label{2.5}
|\lambda| \left\Vert u \right\Vert_{L^2} +\mu\left\Vert u \right\Vert_{H^2}  \leq C \left\Vert f \right\Vert_{L^2},
\end{equation}
where the constant $C>0$ depends only on $\theta, \alpha.$
\end{lemma}
\begin{proof}
Note that $u_2$ satisfies the Dirichlet boundary conditions at $y=0,1$, which is the same as in \cite{Heck}, then the conclusions hold for $u_2$ and we only need to claim that the conclusions hold for $u_1$.

Similar to the arguments in \cite{Heck}, thanks to the Helmholtz projection, we only need to consider the following problem
\begin{equation}\label{2.51}
\left \{
\begin{array}{ll}
(\lambda -\Delta)u_1=f_1 \ \ \mathrm{in} \ \ \Omega,\\
\frac{\alpha}{\mu} u_1(x,0)-\partial_yu_1(x,0)=0,\\
u_1(x,1)=0.
\end{array}
\right.
\end{equation}
Applying the Fourier series, one has
\begin{equation}\label{2.52}
\left \{
\begin{array}{ll}
(\zeta^2 -\partial_y^2)\hat{u}_1=\hat{f}_1, \ 0<y<1,\\
\frac{\alpha}{\mu} \hat{u}_1(0)-\partial_y\hat{u}_1(0)=0,\\
\hat{u}_1(1)=0,
\end{array}
\right.
\end{equation}
where $\zeta=\zeta(k)$ is the unique $\zeta \in \Sigma\left(\frac{\pi-\theta}{2}\right)$ such that $\zeta^2=\lambda+k^2$ and it is easy to see that $\zeta^2 \in \Sigma(\pi-\theta) \subset \mathbb{C}\setminus(-\infty,0]$ (see \cite{Heck} for details).

It follows from the theory of ordinary differential equations that the solution of (\ref{2.52}) can be given by
\begin{equation}\label{2.53}
\hat{u}_1(y)=\int_{0}^{1} G(y,s)\hat{f}_1(s) \mathrm{d}s,
\end{equation}
in which
\begin{equation}\label{2.54}
\begin{aligned}
&G(y,s)=\\
&\frac{\left(\frac{\alpha}{\mu}+\zeta\right)e^{-\zeta(2-s-y)}
+\left(\frac{\alpha}{\mu}-\zeta\right)e^{-\zeta(s+y)}
-\left(\frac{\alpha}{\mu}+\zeta\right)e^{-\zeta|s-y|}
-\left(\frac{\alpha}{\mu}-\zeta\right)e^{-\zeta(2-|s-y|)}}{2\zeta \left[\left(\frac{\alpha}{\mu}+\zeta\right)-\left(\frac{\alpha}{\mu}-\zeta\right)e^{-2\zeta}\right]}
\end{aligned}
\end{equation}
is the Green function of (\ref{2.52}).

It is easy to obtain that
\begin{equation}\label{2.55}
\begin{aligned}
&-\frac{\left(\frac{\alpha}{\mu}+\zeta\right)e^{-\zeta|s-y|}
+\left(\frac{\alpha}{\mu}-\zeta\right)e^{-\zeta(2-|s-y|)}}{2\zeta \left[\left(\frac{\alpha}{\mu}+\zeta\right)-\left(\frac{\alpha}{\mu}-\zeta\right)e^{-2\zeta}\right]}\\
=&
-\frac{\left(\frac{\alpha}{\mu}-\zeta\right)e^{-\zeta(2-s+y)}}
{2\zeta \left[\left(\frac{\alpha}{\mu}+\zeta\right)-\left(\frac{\alpha}{\mu}-\zeta\right)e^{-2\zeta}\right]}\\
& -\frac{\left(\frac{\alpha}{\mu}-\zeta\right)e^{-\zeta(2+s-y)}}{2\zeta \left[\left(\frac{\alpha}{\mu}+\zeta\right)-\left(\frac{\alpha}{\mu}-\zeta\right)e^{-2\zeta}\right]}
-\frac{e^{-\zeta|s-y|}}{2\zeta}\\
:=&G_3+G_4+G_5,
\end{aligned}
\end{equation}
then we have
$$G(y,s)=G_1+G_2+G_3+G_4+G_5,$$
where
$$G_1=\frac{\left(\frac{\alpha}{\mu}+\zeta\right)e^{-\zeta(2-s-y)}}{{2\zeta \left[\left(\frac{\alpha}{\mu}+\zeta\right)-\left(\frac{\alpha}{\mu}-\zeta\right)e^{-2\zeta}\right]}}, \ G_2=\frac{\left(\frac{\alpha}{\mu}-\zeta\right)e^{-\zeta(s+y)}}{{2\zeta \left[\left(\frac{\alpha}{\mu}+\zeta\right)-\left(\frac{\alpha}{\mu}-\zeta\right)e^{-2\zeta}\right]}}.$$

Note that the above Green function (\ref{2.54}) and each term $G_i(i=1,2,3,4,5)$ of the Green function (\ref{2.54}) have the forms which are similar to that of \cite{Heck}, therefore every $G_i(i=1,2,3,4,5)$ can be estimated by the similar arguments as in \cite{Heck}. The rest estimates of this proof can be obtained by the theory of the Fourier multiplier, and we omit it here and refer \cite{Heck} for details.
\end{proof}

It follows from Lemma \ref{lemma2.3} that the Stokes operator $-A$ generates an analytic semigroup $\{e^{-tA}\}$ in $L^2_\sigma$ and $\mathbb{C}\setminus(-\infty,0] \subset \rho(-A)$. In particular, the estimate (\ref{2.5}) holds for $\lambda \in (0,+\infty)$, which implies that $0 \in \rho(-A)$ by some standard arguments, and therefore we can get the classical Stokes estimate
$$\left\Vert (-A)^{-1} f \right\Vert_{H^2}  \leq C \left\Vert f \right\Vert_{L^2}.$$

Define
$$Bu=P\left(u \cdot \nabla v_s+v_s \cdot \nabla u \right), \ \ u \in D(B)=H^1_{*,\sigma}$$
and
$$Lu=(A+B)u, \ \ u \in D(L)=D(A).$$

Recall that the Stokes operator $-A$ generates a $C_0-$semigroup $\{e^{-tA}\}$ in $L^2_\sigma$, which is analytic and bounded in $\Sigma(\theta)$ for $\theta \in (0,\frac{\pi}{2})$. Then for any $f \in D(A), \ \eta>0$, by the interpolation inequality and Poincar$\acute{\mathrm{e}}$'s inequality, we get
$$\left\Vert -B f \right\Vert_{L^2} \leq C \left\Vert u \right\Vert_{H^1} \leq \eta \left\Vert u \right\Vert_{H^2}+C(\eta)\left\Vert u \right\Vert_{L^2},$$
which implies that the operator $-B$ is $(-A)$-bounded and the $(-A)$-bound is 0.
Then the perturbation theory of operator (see \cite{Engel,Kato2} for details) yields that the operator $-L$ generates an analytic semigroup $\{e^{-tL}\}$ in $L^2_\sigma$. Moreover, there exists $\eta_0 >0$ such that for any $f\in L^2_\sigma$ and each $\lambda \in \Sigma(\pi-\theta) \cap \{\lambda \in \mathbb{C}:|\lambda| \geq \eta_0\}, \theta \in(0,\pi)$, the following estimate holds
$$\left\Vert (\lambda+L)^{-1} f \right\Vert_{H^2}  \leq C \left\Vert f \right\Vert_{L^2}.$$
Note that $H^1_{*,\sigma} \cap H^2 \hookrightarrow \hookrightarrow L^2_{\sigma}$, then the operator $(\lambda+L)^{-1}$ is compact in $L^2_{\sigma}$. Hence, $\sigma (-L)$ consists of the isolated eigenvalues of $-L$ and has no accumulation points except infinity.

The following lemma will be used in our analysis.

\begin{lemma}\label{lemma2.2}
For any $f(y)\in H^2(0,1)$ with $f(0)=0$ and $f(1)=0$, there holds
\begin{equation}\label{2.6}
\int_{0}^{1} |f(y)|^2 \mathrm{d}y \leq \int_{0}^{1} |f'(y)|^2 \mathrm{d}y \leq \int_{0}^{1} |f''(y)|^2 \mathrm{d}y.
\end{equation}
\end{lemma}

\begin{proof}
This lemma follows straightforward from integrating by parts, Poincar$\acute{\mathrm{e}}$'s inequality and Young's inequality.
\end{proof}
\begin{remark}
In fact, similar to the proof of the Poincar$\acute{\mathrm{e}}$'s inequality, one can deduce that the Poincar$\acute{\mathrm{e}}$'s inequality holds if $u \in H^1_*$.
\end{remark}

\section{Proof of Theorem \ref{the1.1}: Linear stability}\label{sec3}

In order to analyse the perturbation problem (\ref{1.3})-(\ref{1.4}), we need to study the Stokes operator and perturbed Stokes operator. In fact, we can consider the following abstract Cauchy problem
\begin{equation}\label{3.0}
\left \{
\begin{array}{ll}
\partial_t u+Lu=f(u) & \mathrm{in} \ \Omega, \\
u|_{t=0}=u_0 \ \ & \mathrm{in} \ \Omega,
\end{array}
\right.
\end{equation}
where
$$Lu=P\left(-\mu \Delta u + u \cdot \nabla v_s+v_s \cdot \nabla u \right)$$
is the linear part and $f(u)=P\left(-u\cdot \nabla u\right)$ is the nonlinear term. The linear operator $L$ can be decomposed into the classical Stokes operator $A$ and the perturbed part $B$.

In order to obtain the stability of the Couette flow, we have to show that the spectrum of the operator $-L$ lies on the left side of the complex plane. Then by the standard theory of semigroup, the linear stability is obtained.

Now we are in a position to state the key lemma for the linear stability.

\begin{lemma}\label{lemma3.1}
Under the assumptions of Theorem \ref{the1.1}, there holds
$$m:=\sup \left\{\mathrm{Re} \ \lambda :\lambda \in \sigma (-L)\right\} \leq -C <0,$$
where the constant $C>0$ depends only on $\alpha, \mu, a, b$.
\end{lemma}

\begin{proof}
Since $H^1_{*,\sigma} \cap H^2 \hookrightarrow \hookrightarrow L^2_{\sigma}$, then the operator $(\lambda+L)^{-1}$ is compact in $L^2_{\sigma}$, and therefore $\sigma (-L)$ consists of the isolated eigenvalues of $-L$ and has no accumulation points except infinity.

For a fixed $\theta \in (0,\frac{\pi}{2})$, there exists suitable $r>0$ such that
$$ \overline{\Sigma\left(\frac{\pi}{2}+\theta\right)} \cap \left\{ \lambda \in \mathbb{C} : |\lambda| \geq r\right\} \subset \rho(-L).$$

Note that $\sigma(-L)$ has no accumulation points in $\{\lambda \in \mathbb{C}:|\lambda| \leq r\}$, then we only need to prove that
$$\mathrm{Re} \ \lambda <0, \ \ \ \forall \lambda \in \sigma (-L).$$

Let $\lambda \in \sigma(-L)$ be any eigenvalue of $-L$ and $u\in H^1_{*,\sigma} \cap H^2,u \not\equiv 0$ be the nontrivial eigenvector of $\lambda$, i.e.,
$$(\lambda + L)u=0.$$
The above equation can be rewritten as
$$P(\lambda u-\mu \Delta u+ u \cdot \nabla v_s+v_s \cdot \nabla u)=0.$$
Thanks to Lemma \ref{lemma2.1}, there exists $p\in H^1$ such that
$$\lambda u-\mu \Delta u+ u \cdot \nabla v_s+v_s \cdot \nabla u=-\nabla p.$$

Standard arguments for the elliptic equations guarantee the regularity of $u,p$.
Then $(u,p)$ solves the following problem
\begin{equation}\label{3.1}
\left \{
\begin{array}{lllll}
\lambda u-\mu \Delta u + u \cdot \nabla v_s+v_s \cdot \nabla u +\nabla p=0 \ \ &\mathrm{in} \ \Omega, \\
\nabla \cdot u =0 \ \ &\mathrm{in} \ \Omega, \\
u=0 \ \ &\mathrm{on} \ \Sigma_1,\\
u_2=0 \ \ &\mathrm{on} \ \Sigma_0, \\
\mu \partial_y u_1-\alpha u_1=0 \ \ &\mathrm{on} \ \Sigma_0.
\end{array}
\right.
\end{equation}

The equations in (\ref{3.1}) can be rewritten componentwise as
\begin{equation}\label{3.2}
\left \{
\begin{array}{lllll}
\lambda u_1-\mu \Delta u_1 + \left( \frac{\alpha (a-b)}{\mu +\alpha} y +\frac{\mu a+\alpha b}{\mu +\alpha} \right) \partial_x u_1 +\frac{\alpha (a-b)}{\mu +\alpha}u_2 +\partial_x p=0 &\mathrm{in} \ \Omega, \\
\lambda u_2-\mu \Delta u_2 + \left( \frac{\alpha (a-b)}{\mu +\alpha} y +\frac{\mu a+\alpha b}{\mu +\alpha} \right)  \partial_x u_2  +\partial_y p=0  &\mathrm{in} \ \Omega,\\
\partial_x u_1+\partial_y u_2=0  &\mathrm{in} \ \Omega.
\end{array}
\right.
\end{equation}

In terms of the Fourier series,
$$u(x,y)=\sum_{k \in J} \hat{u}_k(y)e^{ikx}, \ \ p(x,y)=\sum_{k \in J} \hat{p}_k(y) e^{ikx},$$
where $\hat{u}_k,\hat{p}_k$ are smooth on $[0,1]$ and $J$ is some finite subsets of $\mathbb{Z}$, we have
\begin{equation}\label{3.3}
\left \{
\begin{array}{lllll}
\lambda \hat{u}_{1,k}-\mu (\partial_y^2 -k^2)  \hat{u}_{1,k}+ik \left( \frac{\alpha (a-b)}{\mu +\alpha} y +\frac{\mu a+\alpha b}{\mu +\alpha} \right) \hat{u}_{1,k} +\frac{\alpha (a-b)}{\mu +\alpha}\hat{u}_{2,k} +ik \hat{p}_{k}=0, \\
\lambda \hat{u}_{2,k}-\mu (\partial_y^2 -k^2)\hat{u}_{2,k}  +ik\left( \frac{\alpha (a-b)}{\mu +\alpha} y +\frac{\mu a+\alpha b}{\mu +\alpha} \right) \hat{u}_{2,k} +\partial_y \hat{p}_{k}=0,\\
ik\hat{u}_{1,k}+\partial_y \hat{u}_{2,k}=0, \\
\end{array}
\right.
\end{equation}
for $y \in [0,1]$ and $k \in \mathbb{Z}.$ Since $u$ is nontrivial in $\Omega$, then there exists $k \in \mathbb{Z}$ such that $\hat{u}_k \not\equiv 0$. Fixing this $k$ and omitting the subscript $k$ from now, one has
\begin{equation}\label{3.4}
\left \{
\begin{array}{lllll}
\lambda \hat{u}_{1}-\mu (\partial_y^2 -k^2)  \hat{u}_{1}+ik \left( \frac{\alpha (a-b)}{\mu +\alpha} y +\frac{\mu a+\alpha b}{\mu +\alpha} \right) \hat{u}_{1} +\frac{\alpha (a-b)}{\mu +\alpha}\hat{u}_{2} +ik \hat{p}=0, \ \ & 0 <y <1,\\
\lambda \hat{u}_{2}-\mu (\partial_y^2 -k^2)\hat{u}_{2}  +ik \left( \frac{\alpha (a-b)}{\mu +\alpha} y +\frac{\mu a+\alpha b}{\mu +\alpha} \right) \hat{u}_{2} +\partial_y \hat{p}=0,\ \ & 0 <y <1,\\
ik\hat{u}_{1}+\partial_y \hat{u}_{2}=0, \ \  & 0 <y <1, \\
\end{array}
\right.
\end{equation}
which satisfy the following boundary conditions
\begin{equation}\label{3.5}
\left \{
\begin{array}{lllll}
\hat{u}_2(0)=\hat{u}_2(1)=\hat{u}_1 (1)=0,\\
\mu \partial_y \hat{u}_1 (0)-\alpha \hat{u}_1 (0)=0.
\end{array}
\right.
\end{equation}

$\mathbf{Case \ 1:} \ k=0.$

If $k=0,$ then $\partial_y \hat{u}_2=0$ due to (\ref{3.4}), which implies that
$$\hat{u}_2 (y)\equiv \mathrm{constant,\ } y\in  [0,1].$$
Then the boundary conditions yield
$$\hat{u}_2 (y)\equiv 0, \ \ \ y \in [0,1],$$
which implies that
\begin{equation}\label{3.6}
\lambda \hat{u}_{1}-\mu \partial_y^2 \hat{u}_{1}=0, \ \ 0 <y <1.
\end{equation}

Multiplying (\ref{3.6}) by $\overline{\hat{u}_1}$, the complex conjugate of $\hat{u}_1$, and multiplying the conjugate of equation (\ref{3.6}) by $\hat{u}_1$, then integrating over $(0,1)$ and using the boundary conditions, one obtains
\begin{equation}\label{3.7}
\mathrm{Re} \ \lambda \int_{0}^{1} |\hat{u}_1|^2 \mathrm{d}y +\mu \int_{0}^{1} |\partial_y \hat{u}_1|^2 \mathrm{d}y +\alpha |\hat{u}_1 (0)|^2=0.
\end{equation}

If (i) of Theorem \ref{the1.1} holds, that is $\alpha \geq 0,$ then for any $\mu >0$, we have
\begin{equation}\label{3.71}
(\mathrm{Re} \ \lambda +\mu )\int_{0}^{1} |\hat{u}_1|^2 \mathrm{d}y  \leq 0,
\end{equation}
where we have used the Poincar$\acute{\mathrm{e}}$'s inequality.
Therefore
\begin{equation}\label{3.72}
\mathrm{Re} \ \lambda \leq -\mu  <0.
\end{equation}

Now we assume that (ii) of Theorem \ref{the1.1} holds. Simple calculations yield that
\begin{equation}\label{3.8}
\begin{aligned}
\alpha |\hat{u}_1 (0)|^2&=\alpha \int_{0}^{1} \partial_y \left[(y-1)|\hat{u}_1|^2\right] \mathrm{d}y\\
&=\alpha\int_{0}^{1} |\hat{u}_1|^2 \mathrm{d}y+\alpha\int_{0}^{1}2(y-1)  \left[\mathrm{Re} \ \hat{u}_1 \mathrm{Re} \ \partial_y \overline{\hat{u}_1}+\mathrm{Im} \ \hat{u}_1 \mathrm{Im} \ \partial_y \overline{\hat{u}_1}\ \right]
\mathrm{d}y \\
&\geq \alpha\int_{0}^{1}  |\hat{u}_1|^2 \mathrm{d}y-|\alpha|\int_{0}^{1}2  \left|\mathrm{Re} \ \hat{u}_1 \mathrm{Re} \ \partial_y \overline{\hat{u}_1}+\mathrm{Im} \ \hat{u}_1 \mathrm{Im} \ \partial_y \overline{\hat{u}_1}\right| \mathrm{d}y \\
&\geq (\alpha-|\alpha|)\int_{0}^{1}  |\hat{u}_1|^2 \mathrm{d}y-|\alpha| \int_{0}^{1}|\partial_y \hat{u}_1|^2 \mathrm{d}y.
\end{aligned}
\end{equation}

Putting (\ref{3.8}) into (\ref{3.7}) shows that

\begin{equation}\label{3.800}
(\mathrm{Re}\ \lambda +\alpha-|\alpha|) \int_{0}^{1} |\hat{u}_1|^2 \mathrm{d}y +(\mu-|\alpha|) \int_{0}^{1} |\partial_y \hat{u}_1|^2 \mathrm{d}y \leq 0.
\end{equation}

Since $\hat{u}_1(1)=0$ and $\mu-|\alpha|>|\alpha|-\alpha\geq 0$, it follows from Poincar$\acute{\mathrm{e}}$'s inequality and Lemma \ref{lemma2.2} that
\begin{equation}\label{3.8001}
(\mathrm{Re} \ \lambda +\mu+\alpha-2|\alpha|) \int_{0}^{1}  |\hat{u}_1|^2 \mathrm{d}y \leq 0,
\end{equation}
which yields
\begin{equation}\label{3.8002}
\mathrm{Re} \ \lambda \leq 2|\alpha|-\alpha-\mu=-(\mu+3\alpha)<0.
\end{equation}

$\mathbf{Case \ 2:} \ k \neq 0.$

Eliminating $\hat{p}$ in (\ref{3.4}), one has
\begin{equation}\label{3.10}
\mu (\partial_y^2-k^2)^2 \hat{u}_2=\left[ik\left( \frac{\alpha (a-b)}{\mu +\alpha} y +\frac{\mu a+\alpha b}{\mu +\alpha} \right)+\lambda\right](\partial_y^2-k^2) \hat{u}_2, \ \ 0<y<1,
\end{equation}
with the following boundary conditions
\begin{equation}\label{3.11}
\left \{
\begin{array}{lllll}
\hat{u}_2(0)=\hat{u}_2(1)=\partial_y \hat{u}_2(1)=0,\\
\partial_y^2 \hat{u}_2 (0)=\frac{\alpha}{\mu} \partial_y \hat{u}_2 (0).
\end{array}
\right.
\end{equation}

Let $\xi=k\alpha (a-b)$. Multiplying (\ref{3.10}) by $\overline{\hat{u}_2}$, the complex conjugate of $\hat{u}_2$, then integrating over $(0,1)$ and using the boundary conditions (\ref{3.11}), we find that
\begin{equation}\label{3.111}
\begin{aligned}
\mu\left(H_2^2+2k^2 H_1^2+k^4H_0^2\right)=&\frac{i\xi}{\mu+\alpha}\int_{0}^{1} \partial_y \overline{\hat{u}_2} \cdot \hat{u}_2 \mathrm{d}y \\
&-\frac{i\xi}{\mu+\alpha}\left(\int_{0}^{1} y|\partial_y\hat{u}_2|^2\mathrm{d}y+k^2H_0^2\right)\\
&+\frac{ik(\mu a+\alpha b)}{\mu+\alpha}\left(H_1^2+k^2H_0^2\right) -\lambda(\xi) \left(H_1^2+k^2H_0^2\right),
\end{aligned}
\end{equation}
where
$$H_2^2=\int_{0}^{1} |\hat{u}_2|^2\mathrm{d}y +\frac{\alpha}{\mu}|\hat{u}_2(0)|^2,\ H_j^2=\int_{0}^{1} |\partial_y^{j}\hat{u}_2|^2\mathrm{d}y, \ j=0,1.$$

It follows from (\ref{3.111}) that
$$\mathrm{Re}\ \lambda(\xi)=\left(\mathrm{Re}\ \left\{\frac{i\xi}{\mu+\alpha}\int_{0}^{1} \partial_y \overline{\hat{u}_2} \cdot \hat{u}_2 \mathrm{d}y\right\}-\mu\left(H_2^2+2k^2 H_1^2+k^4H_0^2\right)\right)\cdot \left(H_1^2+k^2H_0^2\right)^{-1}.$$

Next, we consider the complex conjugate of the equation (\ref{3.10}):
\begin{equation}\label{3.112}
\mu (\partial_y^2-k^2)^2 \overline{\hat{u}_2}=\left[i\left(\frac{-\xi}{\mu+\alpha}y-\frac{\mu a+\alpha b}{\mu +\alpha}\right)+\overline{\lambda}\right](\partial_y^2-k^2) \overline{\hat{u}_2}, \ \ 0<y<1.
\end{equation}

Multiplying (\ref{3.112}) by $\hat{u}_2$, integrating over $(0,1)$ and using the boundary conditions, similar to (\ref{3.111}), one can get
\begin{equation}\label{3.113}
\begin{aligned}
&\mathrm{Re}\ \overline{\lambda(-\xi)}\\
=&\left(\mathrm{Re}\ \left\{\frac{i\xi}{\mu+\alpha}\int_{0}^{1} \partial_y \overline{\hat{u}_2} \cdot \hat{u}_2 \mathrm{d}y\right\}-\mu\left(H_2^2+2k^2 H_1^2+k^4H_0^2\right)\right)\cdot \left(H_1^2+k^2H_0^2\right)^{-1}\\
=& \mathrm{Re}\ \lambda(\xi).
\end{aligned}
\end{equation}

From these discussions, we can suppose that $\xi=k\alpha (a-b)\geq 0,$ that is, we can always assume that $k>0$ for $\alpha(a-b) \geq 0$ and $k<0$ if $\alpha(a-b) < 0$. Therefore, for simplicity, we rewrite $\xi$ as $\xi=k|\alpha(a-b)| \geq 0, k>0$.

Setting
$$\lambda=-ik\left(\frac{|\alpha(a-b)|}{\mu+\alpha}c+\frac{\mu a+\alpha b}{\mu +\alpha}\right), \ c \in \mathbb{C}, \ R_1=\frac{|\alpha(a-b)|}{\mu(\mu+\alpha)},$$
we obtain the Orr-Sommerfeld boundary value problem
\begin{equation}\label{3.12}
\left \{
\begin{array}{lllll}
(\partial_y^2-k^2)^2 \phi=ikR_1(y-c)(\partial_y^2-k^2) \phi, \ \ 0<y<1,\\
\phi(0)=\phi(1)=\phi'(1)=0,\\
\phi'' (0)=\frac{\alpha}{\mu} \phi' (0),
\end{array}
\right.
\end{equation}
where we have replaced $\hat{u}_2$ by $\phi$ and $\partial_y$ with $'$ for simplicity. Note that $\mathrm{Re}\ \lambda = k\frac{|\alpha(a-b)|}{\mu(\mu+\alpha)}\mathrm{Im}\ c,$ then it suffices to show that the eigenvalue $c \in \mathbb{C}$ of Orr-Sommerfeld problem (\ref{3.12}) satisfies $\mathrm{Im}\ c<0.$

Multiplying (\ref{3.12})$_1$ by $\overline{\phi}$, the complex conjugate of $\phi$, then integrating over $(0,1)$ and using the boundary conditions, one obtains that
\begin{equation}\label{3.13}
\begin{aligned}
\mathrm{Im}\ c =\frac{Q-\overline{Q}-(kR_1)^{-1}\left(I_2^2+2k^2I_1^2+k^4I_0^2\right)}{I_1^2+k^2I_0^2},
\end{aligned}
\end{equation}
where
$$I_2^2=\int_{0}^{1}|\phi''|^2\mathrm{d}y+\frac{\alpha}{\mu}|\phi'(0)|^2, \ I_j^2=\int_{0}^{1}|\phi^{(j)}|^2\mathrm{d}y, \ j=0,1, \ Q=\frac{i}{2}\int_{0}^{1}\phi \overline{\phi'} \mathrm{d}y.$$

By the H$\ddot{\mathrm{o}}$lder's inequality, it holds that
\begin{equation}\label{3.14}
\begin{aligned}
\mathrm{Im}\ c \leq \frac{I_0 I_1-(kR_1)^{-1}\left(I_2^2+2k^2I_1^2+k^4I_0^2\right)}{I_1^2+k^2I_0^2}.
\end{aligned}
\end{equation}

If (i) of Theorem \ref{the1.1} holds, note that $\alpha \geq 0$ and $k>0$, we have
\begin{equation}\label{3.141}
\begin{aligned}
\mathrm{Im}\ c & =\frac{Q-\overline{Q}-(kR_1)^{-1}\left(I_2^2+2k^2I_1^2+k^4I_0^2\right)}{I_1^2+k^2I_0^2} \\
& \leq \frac{Q-\overline{Q}-(kR_1)^{-1}\left(\int_{0}^{1}|\phi''|^2\mathrm{d}y+2k^2I_1^2+k^4I_0^2\right)}{I_1^2+k^2I_0^2}=
:\mathrm{Im}\ \tilde{c}.
\end{aligned}
\end{equation}
Then following the arguments of Romanov in \cite{Romanov}, one can get
$$\mathrm{Im}\ \tilde{c}  <0,$$
then
\begin{equation}\label{3.142}
\begin{aligned}
\mathrm{Im}\ c <0
\end{aligned}
\end{equation}
for any $k>0$, $\alpha \geq 0, a,b \in \mathbb{R}$ and $\mu >0$.

If (ii) of Theorem \ref{the1.1} holds, we need some further estimates as follows.

For $I_j,\ j=0,1,2$, one has
\begin{equation}\label{3.15}
\begin{aligned}
I_2^2&=\int_{0}^{1}|\phi''|^2\mathrm{d}y+\frac{\alpha}{\mu}|\phi'(0)|^2 \\
&=\int_{0}^{1}|\phi''|^2\mathrm{d}y+\frac{\alpha}{\mu}\int_{0}^{1} \left((y-1) |\phi'|^2\right)'\mathrm{d}y \\
& \geq \int_{0}^{1}|\phi''|^2\mathrm{d}y+\frac{\alpha}{\mu}\int_{0}^{1}|\phi'|^2\mathrm{d}y
-\frac{|\alpha|}{\mu}\int_{0}^{1}|\phi'|^2\mathrm{d}y-\frac{|\alpha|}{\mu} \int_{0}^{1}|\phi''|^2\mathrm{d}y \\
&=(1-\frac{|\alpha|}{\mu})\int_{0}^{1}|\phi''|^2\mathrm{d}y+\frac{\alpha-|\alpha|}{\mu}\int_{0}^{1}|\phi'|^2\mathrm{d}y \\
&\geq (1-\frac{2|\alpha|-\alpha}{\mu})\int_{0}^{1}|\phi'|^2\mathrm{d}y=(1-\frac{2|\alpha|-\alpha}{\mu})I_1^2
\end{aligned}
\end{equation}
for $\mu>2|\alpha|-\alpha$, where Lemma \ref{lemma2.2} and Young's inequality have been used.

Similar calculations and the Poincar$\acute{\mathrm{e}}$'s inequality yield that
\begin{equation}\label{3.16}
\begin{aligned}
I_2^2 \geq (1-\frac{2|\alpha|-\alpha}{\mu})I_0^2,
\end{aligned}
\end{equation}
and the classical Poincar$\acute{\mathrm{e}}$'s inequality yields
\begin{equation}\label{3.17}
\begin{aligned}
I_1^2 \geq I_0^2.
\end{aligned}
\end{equation}

Despite (\ref{3.15})--(\ref{3.17}), it seems still difficult to find a useful exact value of the lower bound for
$$\frac{(kR_1)^{-1}\left(I_2^2+2k^2I_1^2+k^4I_0^2\right)}{I_1^2+k^2I_0^2}.$$

To overcome this difficulty, we come out the following analysis.

Let $\delta_0 \in (0,1)$ be given by $2\delta_0^3=1-\delta_0$. Furthermore, for any fixed $\delta \in (\delta_0,1],$ one has
\begin{equation}\label{3.18}
\begin{aligned}
&\frac{I_2^2+2k^2 I_1^2+k^4I_0^2}{I_0I_1} \\
&= \frac{I_2^2}{I_0I_1}+\frac{2k^2}{I_0I_1}\left(\delta I_1^2+(1-\delta) I_1^2+k^2\frac{I_0^2}{2}\right)\\
&= \frac{I_2^2}{I_0I_1} +\frac{2k^2}{I_0I_1}\left[\delta I_1^2+(1-\delta)(I_1-\frac{k(1-\delta)^{-\frac{1}{2}}}{\sqrt{2}}I_0)^2
+\sqrt{2}k(1-\delta)^{\frac{1}{2}}I_0I_1\right]\\
& \geq (1-\frac{2|\alpha|-\alpha}{\mu}) +\frac{2k^2}{I_0I_1}
\max\left\{\sqrt{2}k(1-\delta)^{\frac{1}{2}}I_0I_1,\delta I_1^2\right\}\\
& \geq (1-\frac{2|\alpha|-\alpha}{\mu})+\max\left\{\frac{2k^2}{I_0I_1}\cdot \sqrt{2}k(1-\delta)^{\frac{1}{2}}I_0I_1,\frac{2k^2}{I_0I_1}\cdot \delta I_1^2 \right\}\\
& \geq \max\left\{(1-\frac{2|\alpha|-\alpha}{\mu})+2\sqrt{2}k^3(1-\delta)^{\frac{1}{2}},
(1-\frac{2|\alpha|-\alpha}{\mu})+2k^2\delta\right\}.
\end{aligned}
\end{equation}

For $k \in (0,+\infty)$, define
\begin{align}\label{3.181}
\begin{aligned}
f(k) =&\frac{1}{k}\max\left\{(1-\frac{2|\alpha|-\alpha}{\mu})+2\sqrt{2}k^3(1-\delta)^{\frac{1}{2}},
(1-\frac{2|\alpha|-\alpha}{\mu})+2k^2\delta\right\}\\
=&
\left\{
\begin{array}{ll}
(1-\frac{2|\alpha|-\alpha}{\mu})\frac{1}{k}+2\sqrt{2}(1-\delta)^{\frac{1}{2}}k^2 &, \ k \geq \frac{\sqrt{2}}{2} \delta (1-\delta)^{-\frac{1}{2}}, \\
(1-\frac{2|\alpha|-\alpha}{\mu})\frac{1}{k}+2\delta k&,\ 0< k \leq \frac{\sqrt{2}}{2} \delta (1-\delta)^{-\frac{1}{2}},
\end{array}
\right.
\end{aligned}
\end{align}
and it is easy to see that $f(k) \in C(0,+\infty)$.

For $k \geq \frac{\sqrt{2}}{2} \delta (1-\delta)^{-\frac{1}{2}}$, we have
\begin{align}\label{3.182}
\begin{aligned}
f'(k) =&\frac{\left(2^{\frac{5}{6}}(1-\delta)^{\frac{1}{6}}k-(1-\frac{2|\alpha|-\alpha}{\mu})^{\frac{1}{3}}\right)
}{k^2} \\
&\times \frac{\left(2^{\frac{5}{3}}(1-\delta)^{\frac{1}{3}}k^2
+2^{\frac{5}{6}}(1-\delta)^{\frac{1}{6}}k(1-\frac{2|\alpha|-\alpha}{\mu})^{\frac{1}{3}}
+(1-\frac{2|\alpha|-\alpha}{\mu})^{\frac{2}{3}}\right)}{k^2} > 0,
\end{aligned}
\end{align}

Hence, on $[\frac{\sqrt{2}}{2} \delta (1-\delta)^{-\frac{1}{2}},+\infty)$, it holds that
\begin{align}\label{3.183}
\begin{aligned}
f(k) \geq & f(\frac{\sqrt{2}}{2}\delta (1-\delta)^{-\frac{1}{2}})  \\
=&\sqrt{2}(1-\frac{2|\alpha|-\alpha}{\mu})\delta^{-1}(1-\delta)^{\frac{1}{2}}
+\sqrt{2}\delta^2(1-\delta)^{-\frac{1}{2}} \\
\geq & 2\sqrt{2 \delta}(1-\frac{2|\alpha|-\alpha}{\mu})^{\frac{1}{2}}.
\end{aligned}
\end{align}

If $0< k \leq \frac{\sqrt{2}}{2} \delta (1-\delta)^{-\frac{1}{2}}$, one can get from the average inequality that
$$f(k) \geq 2\sqrt{2 \delta}(1-\frac{2|\alpha|-\alpha}{\mu})^{\frac{1}{2}}.$$

Putting these estimates together leads to
\begin{equation}\label{3.191}
\begin{aligned}
\frac{1}{k}&\cdot \frac{I_2^2+2k^2 I_1^2+k^4I_0^2}{I_0I_1} \\
\geq &
\frac{1}{k}\max\left\{(1-\frac{2|\alpha|-\alpha}{\mu})+2\sqrt{2}k^3(1-\delta)^{\frac{1}{2}},
(1-\frac{2|\alpha|-\alpha}{\mu})+2k^2\delta\right\} \\
\geq & 2\sqrt{2 \delta}(1-\frac{2|\alpha|-\alpha}{\mu})^{\frac{1}{2}}
\end{aligned}
\end{equation}
for $k>0.$

Taking the supremum on the both sides of (\ref{3.191}) on $\delta \in (\delta_0,1]$ gives that
\begin{equation}\label{3.192}
\begin{aligned}
\frac{1}{k}\cdot \frac{I_2^2+2k^2 I_1^2+k^4I_0^2}{I_0I_1} \geq 2\sqrt{2 }(1-\frac{2|\alpha|-\alpha}{\mu})^{\frac{1}{2}}.
\end{aligned}
\end{equation}

Finally, combining (\ref{3.14}), (\ref{3.18}), (\ref{3.191}) and (\ref{3.192}), we obtain
\begin{equation}\label{3.20}
\begin{aligned}
\mathrm{Im}\ c  \leq & \frac{I_0 I_1-(kR_1)^{-1}\left(I_2^2+2k^2I_1^2+k^4I_0^2\right)}{I_1^2+k^2I_0^2} \\
= &\frac{R_1^{-1}I_0I_1}{I_1^2+k^2I_0^2} \left(R_1 -\frac{1}{k}\cdot \frac{I_2^2+2k^2 I_1^2+k^4I_0^2}{I_0I_1}\right)\\
\leq &\frac{R_1^{-1}I_0I_1}{I_1^2+k^2I_0^2} \left(\frac{|\alpha(a-b)|}{\mu(\mu+\alpha)} -2\sqrt{2 }(1-\frac{2|\alpha|-\alpha}{\mu})^{\frac{1}{2}}\right) \\
< & 0
\end{aligned}
\end{equation}
for $\mu>-3\alpha$ and $\frac{|\alpha(a-b)|}{\mu(\mu+\alpha)}\cdot(1+\frac{3\alpha}{\mu})^{-\frac{1}{2}}<2\sqrt{2}$, which completes the proof.
\end{proof}

\section{Proof of Theorem \ref{the1.1}: Nonlinear stability}\label{sec4}

Now we consider the nonlinear problem (\ref{1.3})-(\ref{1.4}). Recall that the nonlinear problem (\ref{1.3})-(\ref{1.4}) can be rewritten as the abstract Cauchy problem
\begin{equation}\label{4.0}
\left \{
\begin{array}{ll}
\partial_t u+Lu=f(u) & \mathrm{in} \ \Omega, \\
u|_{t=0}=u_0 \ \ & \mathrm{in} \ \Omega,
\end{array}
\right.
\end{equation}
where
$$Lu=P\left(-\mu \Delta u + u \cdot \nabla v_s+v_s \cdot \nabla u \right),\ f(u)=P\left(-u\cdot \nabla u\right).$$

To prove Theorem \ref{the1.1}, we need some estimates for fractional powers of operators.

Define
$$A_1=I-P\Delta \ \ \mathrm{with} \ D(A_1)=H^1_{*,\sigma} \cap H^2.$$

Since the operator $A=-P\Delta$ is the generator of an analytic semigroup, then one can define the fractional power of $A_1$. Obviously, the operator $A_1$ is self-adjoint and it is easy to see that the norm $\Vert A_1^{\frac{1}{2}} u \Vert_{L^2}$ is equivalent to $\left\Vert  u \right\Vert_{H^1}$, that is,
\begin{equation}\label{4.000}
\begin{aligned}
\Vert A_1^{\frac{1}{2}} u \Vert_{L^2}\thicksim \left\Vert  u \right\Vert_{H^1}.
\end{aligned}
\end{equation}

The fractional powers of $A_1$ can be estimated as the following lemma.

\begin{lemma}\label{lemma4.1}
There holds
$$ \left\Vert  u \right\Vert_{W^{1,p}} \leq C\left\Vert A_1^{\gamma} u \right\Vert_{L^2}$$
for $u \in D(A_1^{\gamma})$, where the constant $C>0$ depends only on $\gamma,p,$ and $ 1-\frac{1}{p}\leq \gamma<1,p\geq 2.$
\end{lemma}
\begin{proof}
The proof of this lemma is straightforward from Gagiardo-Nirenberg's inequality, H$\ddot{\mathrm{o}}$lder's inequality and Sobolev's inequality, which is similar to the proof of Lemma 5 of \cite{Romanov}. See \cite{Romanov} for details.
\end{proof}

By the arguments similar to Romanov \cite{Romanov}, one can define $A_0:=(sI+L)$ with $D(A_0)=D(A_1)$, where $s=s(\mu,\alpha)>0$ is large enough. For $\gamma \in (0,1)$, define $A_0^{\gamma}$ and the operator $A_0^{\gamma}$ has the equivalent norm
\begin{equation}\label{4.001}
\begin{aligned}
\left\Vert A_0^{\gamma} u \right\Vert_{L^2}\thicksim \left\Vert A_1^{\gamma} u \right\Vert_{L^2}.
\end{aligned}
\end{equation}

Therefore, the Lemma \ref{lemma4.1} holds for $A_0^{\gamma}$, that is
\begin{equation}\label{4.02}
\begin{aligned}
\left\Vert  u \right\Vert_{W^{1,p}} \leq C(\gamma,p)\left\Vert A_0^{\gamma} u \right\Vert_{L^2}, \ u \in D(A_0^{\gamma}),\ 1-\frac{1}{p}\leq \gamma<1,p\geq 2.
\end{aligned}
\end{equation}

Moreover, the following estimate holds (see Romanov \cite{Romanov}):
\begin{equation}\label{4.2}
\begin{aligned}
\left\Vert A_0^{\gamma} e^{-Lt} u \right\Vert_{L^2} \leq C(\mu,\beta,\gamma) t^{-\gamma}e^{-\beta t}\left\Vert u \right\Vert_{L^2}, \ \ \gamma \geq 0, \ t>0, \ \forall \beta \in (0,-m),
\end{aligned}
\end{equation}
where $m$ is defined as in Lemma \ref{lemma3.1}.

Now we are ready to prove Theorem \ref{the1.1}.

{\bf Proof of Theorem \ref{the1.1}.}

By Dehammel's principle, the solution of problem (\ref{1.3})-(\ref{1.4}) is given by
\begin{equation}\label{4.3}
\begin{aligned}
u(t)=e^{-tL}u_0-\int_{0}^{t} e^{-L(t-s)} P(u\cdot \nabla u)(s) \mathrm{d}s.
\end{aligned}
\end{equation}

Define the Picard's sequence
\begin{equation}\label{4.4}
\begin{aligned}
u_n(t)=e^{-tL}u_0- \int_{0}^{t} e^{-L(t-s)} P(u_{n-1} \cdot \nabla u_{n-1})(s) \mathrm{d}s, \ n=1,2,\cdots,
\end{aligned}
\end{equation}
where $u_0 \in D(A_0^{\frac{1}{2}})$.

We define the working space
$$X:=\left\{ u \in D(L): \sup_{t>0} t^{\frac{1}{4}} e^{\beta t} \Vert A_0^{\frac{3}{4}} u(t) \Vert_{L^2} < \infty \right\}$$
with the norm
$$\left\Vert u \right\Vert_{X}=\sup_{t>0} t^{\frac{1}{4}} e^{\beta t} \Vert A_0^{\frac{3}{4}} u(t) \Vert_{L^2}.$$

It is easy to check that $X$ is a Banach space. Next, we only need to show that $\left\Vert u_n \right\Vert_{X}$ is uniformly bounded if $\Vert  A_0^{\frac{1}{2}}u_0 \Vert_{L^2} \leq \varepsilon$ for some small enough $\varepsilon >0$.

It follows from the Sobolev's inequality and (\ref{4.02}) that
\begin{equation}\label{4.5}
\begin{aligned}
\left\Vert P(u\cdot \nabla w) \right\Vert_{L^2} &\leq C \left\Vert u \right\Vert_{L^4}  \left\Vert \nabla w \right\Vert_{L^4}\\
&\leq C  \left\Vert u \right\Vert_{W^{1,\frac{8}{3}}}\left\Vert w \right\Vert_{W^{1,\frac{8}{3}}}\\
&\leq C \Vert A_0^{\frac{3}{4}} u\Vert_{L^2}\Vert A_0^{\frac{3}{4}} w \Vert_{L^2}
\end{aligned}
\end{equation}
for any  $u,w \in D(A_0)$. Then due to (\ref{4.2}), one has
\begin{equation}\label{4.6}
\begin{aligned}
\Vert &A_0^{\frac{3}{4}} u_n(t) \Vert_{L^2} \\
& \leq \Vert A_0^{\frac{3}{4}} e^{-tL}u_0 \Vert_{L^2} +\int_{0}^{t} \Vert A_0^{\frac{3}{4}} e^{-L(t-s)} P(u_{n-1}\cdot \nabla u_{n-1})(s) \Vert_{L^2}\mathrm{d}s \\
&\leq \Vert A_0^{\frac{3}{4}} e^{-tL}u_0 \Vert_{L^2} +C \int_{0}^{t} (t-s)^{-\frac{3}{4}} e^{-\beta(t-s)} \Vert P(u_{n-1}\cdot \nabla u_{n-1})(s) \Vert_{L^2}\mathrm{d}s \\
&\leq Ct^{-\frac{1}{4}} e^{-\beta t} \Vert  A_0^{\frac{1}{2}}u_0 \Vert_{L^2} +C \int_{0}^{t} (t-s)^{-\frac{3}{4}} \Vert A_0^{\frac{3}{4}}u_{n-1} \Vert_{L^2}^2\mathrm{d}s ,
\end{aligned}
\end{equation}
which yields that
\begin{equation}\label{4.7}
\begin{aligned}
\left\Vert u_n \right\Vert_{X} \leq C \Vert  A_0^{\frac{1}{2}}u_0 \Vert_{L^2}+C\left\Vert u_{n-1} \right\Vert_{X}^2.
\end{aligned}
\end{equation}

Then if $\Vert  A_0^{\frac{1}{2}}u_0 \Vert_{L^2} \leq C \Vert u_0 \Vert_{H^1}  \leq \varepsilon$ for some small $\varepsilon >0$, we have
\begin{equation}\label{4.8}
\begin{aligned}
\left\Vert u_n \right\Vert_{X} \leq C \Vert  A_0^{\frac{1}{2}}u_0 \Vert_{L^2} \leq C,
\end{aligned}
\end{equation}
which implies that $\left\Vert u_n \right\Vert_{X}$ is uniformly bounded. Since the embedding $D(A_0^{\frac{3}{4}}) \hookrightarrow D(A_0^{\frac{1}{2}})$ is compact, hence there exists subsequence converges strongly to $u$, which is the global solution of (\ref{1.6}). In addition, it is easy to deduce that $u \in H^1$ from the equivalent norms (\ref{4.000}) and (\ref{4.001}).

Moreover, it follows from the above estimates and (\ref{4.3}) that
\begin{equation}\label{4.9}
\begin{aligned}
\left\Vert A_0^{\gamma}u(t) \right\Vert_{L^2} \leq C t^{\frac{1}{2}-\gamma}e^{-\beta t}\Vert  A_0^{\frac{1}{2}}u_0 \Vert_{L^2} ,\ \ \frac{1}{2} \leq \gamma <1.
\end{aligned}
\end{equation}

Furthermore,
\begin{equation}\label{4.10}
\begin{aligned}
\Vert A_0^{\frac{1}{2}}u(t) \Vert_{L^2} \leq C e^{-\beta t}\Vert  A_0^{\frac{1}{2}}u_0 \Vert_{L^2} ,
\end{aligned}
\end{equation}
which yields that
\begin{equation}\label{4.10a}
\begin{aligned}
\left\Vert u(t) \right\Vert_{H^1} \leq C e^{-\beta t} \left\Vert u_0 \right\Vert_{H^1}.
\end{aligned}
\end{equation}

Therefore Theorem \ref{the1.1} follows.

\section{Proof of Theorem \ref{the1.2}}\label{sec5}

In this section, we will prove Theorem \ref{the1.2}.
Define
$$H^k_{**}:=\{u \in H^k: u \ \mathrm{satisfies \ the \ boundary \ conditions \ in} \ (\ref{1.6})\}$$
and
$$H^k_{**,\sigma}:=\{u \in H^k_{\sigma}: u \ \mathrm{satisfies \ the \ boundary \ conditions \ in} \ (\ref{1.6})\}.$$

One should note that the Lemma \ref{lemma2.3} still holds for the Navier boundary conditions in (\ref{1.6}). More precisely, we have the following lemma.
\begin{lemma}\label{lemma5.01}
Suppose that $\theta\in (0,\frac{\pi}{2})$ and $\lambda \in \Sigma(\frac{\pi}{2}+\theta)$. Then for any $f \in L^2_\sigma$, there exists unique $u \in H^1_{**,\sigma} \cap H^2$ such that
\begin{equation}\label{5lemma1}
(\lambda+A)u=f,
\end{equation}
and the following estimate holds
\begin{equation}\label{5lemma2}
|\lambda| \left\Vert u \right\Vert_{L^2} +\mu\left\Vert u \right\Vert_{H^2}  \leq C \left\Vert f \right\Vert_{L^2},
\end{equation}
where the constant $C>0$ depends only on $\theta, \alpha_l(l=0,1).$
\end{lemma}

\begin{proof}
The proof of this lemma is similar to that of Lemma \ref{lemma2.3}, and we omit it here.
\end{proof}

\begin{lemma}\label{lemma5.1}
Suppose that $u \in H^1_{**,\sigma}$. Then there holds
\begin{equation}\label{5.3}
\begin{aligned}
\left\Vert  \hat{u}_1 \right\Vert_{L^2(0,1)} \leq C_P \left\Vert \partial_y \hat{u}_1 \right\Vert_{L^2(0,1)},
\end{aligned}
\end{equation}
where $C_P>0$ is the best constant so that Poincar$\acute{\mathrm{e}}$ inequality for $f(y)\in H^1(0,1)$ with $\int_0^1f(y) \mathrm{d}y=0$ holds, $\hat{u}_1$ is defined as before.
\end{lemma}

\begin{proof}
Since $\nabla \cdot u=0,$ thus
$$ik \hat{u}_1+\partial_y \hat{u}_2=0.$$

Note that $u \in H^1_{**,\sigma}$, then $ \hat{u}_2(0)=\hat{u}_2(1)=0$ and
$$\int_{0}^{1} ik \hat{u}_1 \mathrm{d}y=-\int_{0}^{1} \partial_y \hat{u}_2 \mathrm{d}y=0,$$
therefore
$$\int_{0}^{1} \hat{u}_1 \mathrm{d}y=0,$$
which implies that the classical Poincar$\acute{\mathrm{e}}$'s inequality holds for $\hat{u}_1$.
\end{proof}

Next, we give an estimate for the $\sigma(-L)$.
\begin{lemma}\label{lemma5.2}
Under the assumptions of Theorem \ref{the1.2}, there holds
$$\sup \left\{\mathrm{Re} \ \lambda :\lambda \in \sigma (-L)\right\} \leq -{\tilde C}<0,$$
where the constant ${\tilde C}>0$ depends only on $\mu, \alpha_l(l=0,1), a, b, \Omega, C_P.$
\end{lemma}

\begin{proof}
Consider the problem
\begin{equation}\label{5.4}
\left \{
\begin{array}{lllll}
\lambda u-\mu \Delta u + u \cdot \nabla v_s+v_s \cdot \nabla u +\nabla p=0 \ \ &\mathrm{in} \ \Omega, \\
\nabla \cdot u =0 \ \ &\mathrm{in} \ \Omega, \\
u_2=0 \ \ &\mathrm{in} \ \Sigma, \\
\mu \partial_y u_1+\alpha_1 u_1=0 \ \ &\mathrm{on} \ \Sigma_1,\\
\mu \partial_y u_1-\alpha_0 u_1=0 \ \ &\mathrm{on} \ \Sigma_0.
\end{array}
\right.
\end{equation}

Similar to Lemma \ref{lemma3.1}, the Lemma \ref{lemma5.01} and the Fourier series give that
\begin{equation}\label{5.5}
\left \{
\begin{array}{lllll}
\lambda \hat{u}_{1}-\mu (\partial_y^2 -k^2)  \hat{u}_{1}
\\
\quad  +ik \left(\frac{\alpha_0 \alpha_1(a-b)}{\mu(\alpha_0+\alpha_1) +\alpha_0\alpha_1} y +\frac{\mu(\alpha_1 a+\alpha_0 b)+\alpha_0\alpha_1 b}{\mu(\alpha_0 +\alpha_1)+\alpha_0\alpha_1}\right)\hat{u}_{1}  +\frac{\alpha_0\alpha_1(a-b)}{\mu(\alpha_0+\alpha_1) +\alpha_0\alpha_1}\hat{u}_{2} +ik \hat{p}=0, \\
\lambda \hat{u}_{2}-\mu (\partial_y^2 -k^2)\hat{u}_{2}  +ik\left(\frac{\alpha_0\alpha_1(a-b)}{\mu(\alpha_0+\alpha_1) +\alpha_0\alpha_1} y +\frac{\mu(\alpha_1 a+\alpha_0 b)+\alpha_0\alpha_1 b}{\mu(\alpha_0 +\alpha_1)+\alpha_0\alpha_1}\right) \hat{u}_{2} +\partial_y \hat{p}=0,\\
ik\hat{u}_{1}+\partial_y \hat{u}_{2}=0,  \\
\end{array}
\right.
\end{equation}
with the following boundary conditions
\begin{equation}\label{5.61}
\left \{
\begin{array}{lllll}
\hat{u}_2(0)=\hat{u}_2(1)=0,\\
\mu \partial_y \hat{u}_1 (1)+\alpha_1 \hat{u}_1 (1)=0,\\
\mu \partial_y \hat{u}_1 (0)-\alpha_0 \hat{u}_1 (0)=0.
\end{array}
\right.
\end{equation}

There are two cases to be considered.

$\mathbf{Case \ 1:} \ k=0.$

If $k=0$, then $\hat{u}_{2} \equiv 0$ at $[0,1].$ Therefore
\begin{equation}\label{5.71}
\begin{aligned}
\lambda \hat{u}_{1}-\mu \partial_y^2 \hat{u}_{1}=0.
\end{aligned}
\end{equation}
Multiplying (\ref{5.71}) by $\overline{\hat{u}_1}$, integrating over $(0,1)$ and using the boundary conditions (\ref{5.61}), one gets
\begin{equation}\label{5.81}
\begin{aligned}
\mathrm{Re} \ \lambda \int_{0}^{1} |\hat{u}_1|^2 \mathrm{d}y+\mu \int_{0}^{1} |\partial_y\hat{u}_1|^2 \mathrm{d}y
+\sum\limits_{l=0}^{1} \alpha_l |\hat{u}_1(l)|^2=0.
\end{aligned}
\end{equation}

If the condition (iii) of Theorem \ref{the1.2} holds, that is, $\alpha_l \geq 0 (l=0,1)$, then one has
\begin{equation}\label{5.811}
\begin{aligned}
\left(\mathrm{Re}\ \lambda+\frac{1}{C_P} \mu \right)  \int_{0}^{1} |\hat{u}_1|^2 \mathrm{d}y \leq 0,
\end{aligned}
\end{equation}
where we have used the Lemma \ref{lemma5.1}. Therefore
\begin{equation}\label{5.812}
\begin{aligned}
\mathrm{Re}\ \lambda \leq -\frac{1}{C_P} \mu <0,
\end{aligned}
\end{equation}
in which $C_P>0$ is the best Poincar$\acute{\mathrm{e}}$ constant in Lemma \ref{lemma5.1}.

Now we suppose that the condition (iv) of Theorem \ref{the1.2} holds.
For the boundary terms, it follows from some simple calculations that
\begin{equation}\label{5.9}
\begin{aligned}
&\sum\limits_{l=0}^{1} \alpha_l |\hat{u}_1(l)|^2 \\
=&\int_{0}^{1} \partial_y \left[((\alpha_1+\alpha_0)y-\alpha_0)|\hat{u}_1|^2\right] \mathrm{d}y\\
=&(\alpha_0+\alpha_1) \int_{0}^{1} |\hat{u}_1|^2 \mathrm{d}y \\
& +\int_{0}^{1}2[((\alpha_0+\alpha_1)y-\alpha_0)]  \left[\mathrm{Re} \ \hat{u}_1 \mathrm{Re} \ \partial_y \overline{\hat{u}_1}+\mathrm{Im} \ \hat{u}_1 \mathrm{Im} \ \partial_y \overline{\hat{u}_1}\ \right]
\mathrm{d}y \\
\geq &(\alpha_1+\alpha_0)\int_{0}^{1}  |\hat{u}_1|^2 \mathrm{d}y-\max\limits_{l=0,1}\{|\alpha_l|\}\int_{0}^{1}2  \left|\mathrm{Re} \ \hat{u}_1 \mathrm{Re} \ \partial_y \overline{\hat{u}_1}+\mathrm{Im} \ \hat{u}_1 \mathrm{Im} \ \partial_y \overline{\hat{u}_1}\right|  \mathrm{d}y \\
\geq& (\alpha_1+\alpha_0-\max\limits_{l=0,1}\{|\alpha_l|\})\int_{0}^{1}  |\hat{u}_1|^2 \mathrm{d}y-\max\limits_{l=0,1}\{|\alpha_l|\} \int_{0}^{1}|\partial_y \hat{u}_1|^2 \mathrm{d}y.
\end{aligned}
\end{equation}
Putting the above estimates into (\ref{5.81}) and using Lemma \ref{lemma5.1} yield that
\begin{equation}\label{5.9a}
\begin{aligned}
\left(\mathrm{Re}\ \lambda +(\alpha_0+\alpha_1-\max\limits_{l=0,1}\{|\alpha_l|\})+\frac{1}{C_P}(\mu-\max\limits_{l=0,1}\{|\alpha_l|\})\right)\int_{0}^{1}  |\hat{u}_1|^2 \mathrm{d}y \leq 0,
\end{aligned}
\end{equation}
which implies that
\begin{equation}\label{5.9b}
\begin{aligned}
\mathrm{Re} \ \lambda \leq \left(\frac{1}{C_P}+1\right)\max\limits_{l=0,1}\{|\alpha_l|\}-(\alpha_0+\alpha_1)-\frac{1}{C_P}\mu :=&-{\tilde C} <0
\end{aligned}
\end{equation}
for $\mu>(1+C_P)\max\limits_{l=0,1}\{|\alpha_l|\}-C_P(\alpha_0+\alpha_1)$, where $C_P>0$ is the best Poincar$\acute{\mathrm{e}}$ constant in Lemma \ref{lemma5.1}.

$\mathbf{Case \ 2:} \ k \neq  0.$

By eliminating $\hat{p}$, one has
\begin{equation}\label{5.00}
\begin{aligned}
\mu (\partial_y^2-k^2)^2& \hat{u}_2\\
=&\left[ik\left(\frac{\alpha_0\alpha_1(a-b)}{\mu(\alpha_0+\alpha_1) +\alpha_0\alpha_1} y +\frac{\mu(\alpha_1 a+\alpha_0 b)+\alpha_0\alpha_1 b}{\mu(\alpha_0 +\alpha_1)+\alpha_0\alpha_1}\right)+\lambda\right](\partial_y^2-k^2) \hat{u}_2
\end{aligned}
\end{equation}
for $0<y<1.$

Similar to Lemma \ref{lemma3.1}, setting
$$\lambda=-ik\left(\left|\frac{\alpha_0\alpha_1(a-b)}{\mu(\alpha_0+\alpha_1) +\alpha_0\alpha_1}\right| c +\frac{\mu(\alpha_1 a+\alpha_0 b)+\alpha_0\alpha_1 b}{\mu(\alpha_0 +\alpha_1)+\alpha_0\alpha_1}\right), \ c \in \mathbb{C}$$
and
$$R_2:=\left|\frac{\alpha_0\alpha_1(a-b)}{\mu\left(\mu(\alpha_0+\alpha_1) +\alpha_0\alpha_1\right)}\right|,$$
one can get
\begin{equation}\label{5.10}
\left \{
\begin{array}{lllll}
(\partial_y^2-k^2)^2 \phi=ikR_2(y-c)(\partial_y^2-k^2) \phi, \ \ 0<y<1,\\
\phi(0)=\phi(1)=0,\\
\phi'' (0)=\frac{\alpha_0}{\mu} \phi' (0),\\
\phi'' (1)=-\frac{\alpha_1}{\mu} \phi' (1),
\end{array}
\right.
\end{equation}
where we have replaced $\hat{u}_2$ by $\phi$ and $\partial_y$ with $'$ for simplicity.

Multiplying (\ref{5.10})$_1$ by $\overline{\phi}$, the complex conjugate of $\phi$, then integrating over $(0,1)$ and using the boundary conditions, yield that
\begin{equation}\label{5.11}
\begin{aligned}
\mathrm{Im}\ c =\frac{Q-\overline{Q}-(kR_2)^{-1}\left(I_2^2+2k^2I_1^2+k^4I_0^2\right)}{I_1^2+k^2I_0^2},
\end{aligned}
\end{equation}
where
$$I_2^2=\int_{0}^{1}|\phi''|^2\mathrm{d}y+\sum_{l=0}^{1}\frac{\alpha_l}{\mu}|\phi'(l)|^2, \ I_j^2=\int_{0}^{1}|\phi^{(j)}|^2\mathrm{d}y, \ j=0,1, \ Q=\frac{i}{2}\int_{0}^{1}\phi \overline{\phi'} \mathrm{d}y.$$

It follows from the H$\ddot{\mathrm{o}}$lder's inequality that
\begin{equation}\label{5.12}
\begin{aligned}
\mathrm{Im}\ c \leq \frac{I_0 I_1-(kR_2)^{-1}\left(I_2^2+2k^2I_1^2+k^4I_0^2\right)}{I_1^2+k^2I_0^2}.
\end{aligned}
\end{equation}

If (iii) of Theorem \ref{the1.2} holds, note that $\alpha_l \geq 0(l=0,1)$ and $k>0$, we have
\begin{equation}\label{5.141}
\begin{aligned}
\mathrm{Im}\ c & =\frac{Q-\overline{Q}-(kR_2)^{-1}\left(I_2^2+2k^2I_1^2+k^4I_0^2\right)}{I_1^2+k^2I_0^2} \\
& \leq \frac{Q-\overline{Q}-(kR_2)^{-1}\left(\int_{0}^{1}|\phi''|^2\mathrm{d}y+2k^2I_1^2+k^4I_0^2\right)}{I_1^2+k^2I_0^2}=
:\mathrm{Im}\ \breve{c}.
\end{aligned}
\end{equation}
The arguments of Romanov in \cite{Romanov} give that
$$\mathrm{Im}\ \breve{c}  <0,$$
therefore one has
\begin{equation}\label{5.142}
\begin{aligned}
\mathrm{Im}\ c <0
\end{aligned}
\end{equation}
for any $k>0$, $\alpha_l \geq 0(l=0,1), a,b \in \mathbb{R}$ and $\mu >0$.

Let us suppose that (iv) of Theorem \ref{the1.2} holds.

Now we estimate $I_j,j=0,1,2$.
For $I_2$ and $I_1$, it holds that
\begin{equation}\label{5.13}
\begin{aligned}
I_2^2=&\int_{0}^{1}|\phi''|^2\mathrm{d}y+\sum_{l=0}^{1}\frac{\alpha_l}{\mu}|\phi'(l)|^2 \\
=&\int_{0}^{1}|\phi''|^2\mathrm{d}y+\frac{1}{\mu}\int_{0}^{1} \left(((\alpha_0+\alpha_1)y-\alpha_0)|\phi'|^2\right)'\mathrm{d}y \\
 \geq & \int_{0}^{1}|\phi''|^2\mathrm{d}y+\frac{\alpha_0+\alpha_1}{\mu}\int_{0}^{1}|\phi'|^2\mathrm{d}y\\
&-\frac{\max\limits_{l=0,1}|\alpha_l|}{\mu}\int_{0}^{1}|\phi'|^2\mathrm{d}y-\frac{\max\limits_{l=0,1}|\alpha_l|}{\mu} \int_{0}^{1}|\phi''|^2\mathrm{d}y \\
=&\left(1-\frac{\max\limits_{l=0,1}|\alpha_l|}{\mu}\right)\int_{0}^{1}|\phi''|^2\mathrm{d}y
+\frac{\alpha_0+\alpha_1-\max\limits_{l=0,1}|\alpha_l|}{\mu}\int_{0}^{1}|\phi'|^2\mathrm{d}y \\
\geq & \left(1-\frac{2\max\limits_{l=0,1}|\alpha_l|-\alpha_0-\alpha_1}{\mu}\right)\int_{0}^{1}|\phi'|^2\mathrm{d}y
=\left(1-\frac{2\max\limits_{l=0,1}|\alpha_l|-\alpha_0-\alpha_1}{\mu}\right)I_1^2.
\end{aligned}
\end{equation}
where Lemma \ref{lemma2.2} has been used again.

Similarly, by the Poincar$\acute{\mathrm{e}}$'s inequality, one obtains from Lemma \ref{lemma2.2} that
\begin{equation}\label{5.14}
\begin{aligned}
I_1^2 \geq I_0^2
\end{aligned}
\end{equation}
and then
\begin{equation}\label{5.15}
\begin{aligned}
I_2^2 \geq \left(1-\frac{2\max\limits_{l=0,1}|\alpha_l|-\alpha_0-\alpha_1}{\mu}\right)I_0^2.
\end{aligned}
\end{equation}

Moreover, for any fixed $\delta \in (\delta_0,1],$ one has
\begin{equation}\label{5.6}
\begin{aligned}
&\frac{I_2^2+2k^2 I_1^2+k^4I_0^2}{I_0I_1} \\
&= \frac{I_2^2}{I_0I_1}+\frac{2k^2}{I_0I_1}\left(\delta I_1^2+(1-\delta) I_1^2+k^2\frac{I_0^2}{2}\right)\\
&= \frac{I_2^2}{I_0I_1} +\frac{2k^2}{I_0I_1}\left[\delta I_1^2+(1-\delta)(I_1-\frac{k(1-\delta)^{-\frac{1}{2}}}{\sqrt{2}}I_0)^2
+\sqrt{2}k(1-\delta)^{\frac{1}{2}}I_0I_1\right]\\
& \geq \left(1-\frac{2\max\limits_{l=0,1}|\alpha_l|-\alpha_0-\alpha_1}{\mu}\right) +\frac{2k^2}{I_0I_1}
\max\left\{\sqrt{2}k(1-\delta)^{\frac{1}{2}}I_0I_1,\delta I_1^2\right\}\\
& \geq h+\max\left\{\frac{2k^2}{I_0I_1}\cdot \sqrt{2}k(1-\delta)^{\frac{1}{2}}I_0I_1,\frac{2k^2}{I_0I_1}\cdot \delta I_1^2 \right\}\\
& \geq \max\left\{h+2\sqrt{2}k^3(1-\delta)^{\frac{1}{2}},
h+2k^2\delta\right\}.
\end{aligned}
\end{equation}
where $\delta_0 \in (0,1)$ is given by $2\delta_0^3=1-\delta_0$ and
$$h=\left(1-\frac{2\max\limits_{l=0,1}|\alpha_l|-\alpha_0-\alpha_1}{\mu}\right).$$

For $k \in (0,+\infty)$, define
\begin{align}\label{5.7}
\begin{aligned}
g(k) =&\frac{1}{k}\max\left\{(1-\frac{2|\alpha|-\alpha}{\mu})+2\sqrt{2}k^3(1-\delta)^{\frac{1}{2}},
(1-\frac{2|\alpha|-\alpha}{\mu})+2k^2\delta\right\}\\
=&
\left\{
\begin{array}{ll}
\frac{h}{k}+2\sqrt{2}(1-\delta)^{\frac{1}{2}}k^2 &, \ k \geq \frac{\sqrt{2}}{2} \delta (1-\delta)^{-\frac{1}{2}}, \\
\frac{h}{k}+2\delta k&,\ 0< k \leq \frac{\sqrt{2}}{2} \delta (1-\delta)^{-\frac{1}{2}},
\end{array}
\right.
\end{aligned}
\end{align}
and it is easy to see that $g(k) \in C(0,+\infty)$.

Similar arguments as in Lemma \ref{lemma3.1} give
\begin{equation}\label{5.8}
\begin{aligned}
\frac{1}{k} \frac{I_2^2+2k^2 I_1^2+k^4I_0^2}{I_0I_1} \geq 2\sqrt{2}\left(1-\frac{2\max\limits_{l=0,1}|\alpha_l|-\alpha_0-\alpha_1}{\mu}\right)^{\frac{1}{2}}.
\end{aligned}
\end{equation}
Furthermore, one has
\begin{equation}\label{5.8a}
\begin{aligned}
\mathrm{Im}\ c  \leq & \frac{I_0 I_1-(kR_2)^{-1}\left(I_2^2+2k^2I_1^2+k^4I_0^2\right)}{I_1^2+k^2I_0^2} \\
= &\frac{R_2^{-1}I_0I_1}{I_1^2+k^2I_0^2} \left(R_2-\frac{1}{k}\cdot \frac{I_2^2+2k^2 I_1^2+k^4I_0^2}{I_0I_1}\right)\\
\leq &\frac{R_2^{-1}I_0I_1}{I_1^2+k^2I_0^2} \left[\left|\frac{\alpha_0\alpha_1(a-b)}{\mu\left(\mu(\alpha_0+\alpha_1) +\alpha_0\alpha_1\right)}\right| -2\sqrt{2}\left(1-\frac{2\max\limits_{l=0,1}|\alpha_l|
-\alpha_0-\alpha_1}{\mu}\right)^{\frac{1}{2}}\right] \\
:=&-{\tilde C}<0
\end{aligned}
\end{equation}
if
$$\mu>2\max\limits_{l=0,1}\{|\alpha_l|\}-(\alpha_0+\alpha_1)$$
and
$$\left|\frac{\alpha_0\alpha_1(a-b)}{\mu\left(\mu(\alpha_0+\alpha_1) +\alpha_0\alpha_1\right)}\right|\cdot \left(1-\frac{2\max\limits_{l=0,1}|\alpha_l|-\alpha_0-\alpha_1}{\mu}\right)^{-\frac{1}{2}} <2\sqrt{2}.$$

This completes the proof.
\end{proof}

Lemma \ref{lemma5.2} implies the linear stability in Theorem \ref{the1.2}. It is turn to prove the nonlinear stability.

{\bf Proof of Theorem \ref{the1.2}.}

With linear stability obtained by Lemma \ref{lemma5.2} at hand, one can prove the nonlinear stability by using similar arguments as in Section \ref{sec4}, and therefore the details are omitted here.


{\bf Acknowledgements.}
The authors would like to thank Professor Zhouping Xin, Professor Yan Guo and Professor Xiaoping Wang for their valuable discussions and suggestions. Ding's research is supported by the National Natural Science Foundation of China (No.11371152, No.11571117, No.11871005 and No.11771155) and Guangdong Provincial Natural Science Foundation (No.2017A030313003).


\end{document}